\magnification=1200
\voffset=-1.1truept
\vsize=22.0truecm
\hoffset=.5truecm
\hsize=14.8truecm

\parindent=16pt
\baselineskip=17truept

\font\cmrnine= cmr9
\font\TitFont = cmr17

\def\qed{\vbox to 7.0pt{\hsize=4.0pt \hrule height 7.0pt width 4.0pt }}
\def\re{\Re \kern-1.4truept e}

\nopagenumbers
\def\makeheadline{\lineskip=20truept \line{\the\headline}}
\def\intestazio{\ifodd\pageno\rightheadline \else\leftheadline\fi}
\def\rightheadline{\cmrnine \hfil\folio}
\def\leftheadline{\cmrnine \folio \hfil }
\headline={\ifnum\pageno=1{\hss}\else\intestazio\fi}

\def\IR{I\!\!R}
\def\IC{I\!\!\!\!C}

\centerline{\TitFont Sufficient Conditions for Conservativity}
\centerline{\TitFont of Minimal Quantum Dynamical Semigroups}
\vskip 2.0truecm
$$\vbox { \settabs  2 \columns
\+ A.M. Chebotarev$^{\dagger}$  & F. Fagnola$^{\ddagger}$   \cr
\+  {\cmrnine Moscow State University  } &
{\cmrnine Universit\`a di Genova } \cr
\+ {\cmrnine Quantum Statistics Department}  &
{\cmrnine Dipartimento di  Matematica} \cr
\+ {\cmrnine  119899 Moscow - RUSSIA }
& {\cmrnine  Via Dodecaneso, 35}\cr
\+{\cmrnine }
& {\cmrnine I - 16146 Genova - ITALY}\cr}
$$

\null\vskip 11truecm
\line{\qquad AMS Classification: 81S25  47N50   46L60\hfil}
\footnote{$ $}{
(${\dagger}$) Supported  by GNAFA-CNR Italy.\quad
(${\ddagger}$) Member of GNAFA-CNR Italy. }

\vfill\eject

\noindent{\bf Abstract}. The conservativity of a minimal
quantum dynamical semigroup is proved whenever there exists
a ``reference'' subharmonic operator bounded from
below by the dissipative part of the infinitesimal generator.
We discuss applications of this criteria
in mathematical physics and quantum probability.

\vskip 0.8truecm
\centerline{\bf 1. Introduction}
\noindent A quantum dynamical semigroup (q.d.s.) ${\cal T}=
({\cal T}_t)_{t\ge 0}$ on ${\cal B}(h)$, the Banach space of bounded
operators on a Hilbert space $h$, is a w$^*$-continuous semigroup of
completely positive linear maps on ${\cal B}(h)$. Here $I$ denotes the
identity operator. A q.d.s. is {\sl conservative} (or identity
preserving,or Markovian) if ${\cal T}_t(I) = I$.

Q.d.s. arise in the study of irreversible evolutions in quantum
mechanics (see [2], [3], [4] and the references therein) and
as a quantum analogue of classical Markovian semigroups in quantum
probability (see [20], [21]).

In rather general cases the infinitesimal generator ${\cal L}$
can be written (formally) as
$$
{\cal L}(X) = i\left[H,X\right]
-{1\over 2}XM + \sum_{\ell=1}^\infty L_\ell^*XL_\ell
-{1\over 2}MX, \eqno(1.1)
$$
where $M=\sum_{\ell=1}^\infty L_\ell^*L_\ell$ and $H$ is a symmetric
operator satisfying some conditions that will be made precise later.
However, even if ${\cal L}(I)=0$ for unbounded generators
(1.1), the q.d.s. with the formal generator (1.1) may not be unique
and conservative (see examples in [6], [7], [12]).

The study of conservativity conditions is important in quantum
probability because they play a key role in the proof of
uniqueness and unitarity of solutions of an Hudson--Parthasarathy
quantum stochastic differential equation (see e.g. [7], [8],
[14], [15], [20], [21]). Moreover they allow to deduce
regularity conditions for trajectories of classical Markov processes
(see, for example, [10], [17]) from an operator-theoretic
approach.

On the other hand, when a q.d.s. is conservative, the predual semigroup
(see [4], [12]) is trace preserving.
A conservative irreversible Markov evolution on a von Neumann operator
algebra can be considered as an analog of an isometric evolution in
Hilbert space, and the infinitesimal generators of conservative
q.d.s. play the same role in our approach as essentially self-adjoint
operators generating unitary groups in Hilbert spaces.

Necessary and sufficient conditions for conservativity of a q.d.s. were
obtained in [7],[8]. Some of these conditions, however, are excessive
and difficult to check in practically interesting examples.
A simplified form of sufficient conditions was described in [9] and
improved essentially in our previous paper [11] in view of
applications to quantum stochastic calculus (see, for example, [15]).
These conditions can be written formally as follows
$$
i[H,M] \le b M,\qquad\qquad {\cal L}(M) \le b M
\eqno(1.2)
$$
where $b$ is a constant.

These results are improved here in several aspects. In fact:
\item{1.} The second inequality (1.2) is assumed only for some
self-adjoint operator $C$ bounded from below by $M$,
\item{2.} The operators $L_\ell$ need not to be closed or even closable,
as well as the form represented by the operator $-M-2iH$
(see Example 5.3),

Moreover the previous proof based on technical inequalities for
contractive completely positive maps is replaced by a simple one
based on a priori bounds for the resolvent of the minimal q.d.s. (see
Theorem 3.1). The operator $C$ in the new sufficient condition can be
considered as a ``generalized'' subharmonic operator for the q.d.s.
${\cal T}$.

In Section 5 we apply our conditions to three q.d.s. The first one
arises in a phenomenological model for a physical evolution (see [3],
[4]); conservativity of this q.d.s. could not be proved by
direct applications of the techniques developed in [9], [11].
The second one is the problem of constructing a quantum extension of the
Brownian motion on $[0,+\infty[$ with partial reflection at $0$ where
a non-closable operator $L_\ell$ appears. In this case the
appropriate operator $C$ turns out to be a singular perturbation of the
second derivative on $[0,+\infty[$ by a delta function of the form
studied in [1]. The third is the simplest example in which the
quadratic form represented by the operator $-M-2iH$ is not closed
(see [18]).

\vskip 0.6truecm
\centerline{\bf 2. The minimal quantum dynamical semigroup}

\noindent Let $h$ be a complex separable Hilbert space endowed with
the scalar product $\left\langle\cdot,\cdot\right\rangle$ and norm
$\left\Vert\cdot\right\Vert$ and let ${\cal B}(h)$ be the Banach
space of bounded operators in $h$. The uniform norm in ${\cal B}(h)$
will be denoted by $\left\Vert\cdot\right\Vert_\infty$ and
the identity operator in $h$ will be denoted by $I$. We shall
denote by $D(G)$ the domain of an operator $G$ in $h$.

\proclaim{Definition {2.1}}. A {\rm quantum dynamical semigroup}
({\rm q.d.s.} in the sequel) in ${\cal B}(h)$ is a family
${\cal T}=\left({\cal T}_t\right)_{t\ge 0}$ of operators in
${\cal B}(h)$ with the following properties:
\item{i)} ${\cal T}_0(X) = X$, for all $X\in{\cal B}(h)$,
\item{ii)} ${\cal T}_{t+s}(X)={\cal T}_t\left({\cal T}_s(X)\right)$, for
all $s,t\ge 0$ and all $X\in{\cal B}(h)$,
\item{iii)} ${\cal T}_t(I)\le I$, for all $t\ge 0$,
\item{iv)} (complete positivity) for all $t\ge 0$, for all integer $n$
and all finite sequences $(X_j)_{j=1}^n$, $(Y_l)_{l=1}^n$ of elements
of ${\cal B}(h)$ we have
$$
\sum_{j,l=1}^n Y_l^*{\cal T}_t(X_l^*X_j)Y_j\ge 0,
$$
\item{v)} (normality) for every sequence $(X_n)_{n\ge 0}$ of elements of
${\cal B}(h)$  converging weakly to an element $X$ of ${\cal B}(h)$
the sequence $({\cal T}_t(X_n))_{n\ge 0}$ converges weakly to
${\cal T}_t(X)$ for all  $t\ge 0$,
\item{vi)} (ultraweak continuity) for all trace class operator $\rho$
in $h$ and all  $X\in{\cal B}(h)$ we have
$$
\lim_{t\to 0^+}{\rm Tr}(\rho{\cal T}_t(X)) =
{\rm Tr}(\rho X).
$$

We recall that:

a) as a consequence of properties iii, iv, (see e.g. [13]),
for all $X\in{\cal B}(h)$ and all $t\ge 0$, we have the inequality
$$
\left\Vert{\cal T}_t(X)\right\Vert_\infty
\le \left\Vert X\right\Vert_\infty.\eqno(2.1)
$$
Thus, for all $t\ge 0$, ${\cal T}_t$ is continuous also for the norm
$\left\Vert\cdot\right\Vert_\infty$.

b) as a consequence of properties iv, vi, for all $X\in {\cal B}(h)$,
the map $t\to {\cal T}_t(X)$ is strongly continuous.

\proclaim{Definition {2.2}}. A q.d.s. ${\cal T}$ is said to be
{\rm conservative} if ${\cal T}_t(I) = I$ for all $t\ge 0$.

The bounded infinitesimal generator of a norm continuous q.d.s.
was characterized by Gorini, Lindblad, Kossakowski and Sudarshan
(see e.g. [21] Th.~30.12 p.~267). The characterization problem for
arbitrary q.d.s. is still open  (see e.g. [16]).

A very large class of q.d.s. with unbounded generators was constructed
by Davies [12] (see also [23]) by considering operators $G$, $L_\ell$
($\ell=1,2,\dots $) satisfying the following technical assumption:

\medskip

{\bf A} -  {\sl The operator $G$ is the infinitesimal generator of a
strongly continuous contraction semigroup $P=(P(t))_{t\ge 0}$ in $h$.
The domain of the operators $(L_\ell)_{\ell=1}^\infty$ contains
the domain $D(G)$ of the operator $G$.
For all $u,v\in D(G)$, we have
$$
\left\langle v,Gu\right\rangle +
\left\langle Gv,u\right\rangle +
\sum_{\ell=1}^\infty
\left\langle L_\ell v,L_\ell u\right\rangle = 0. \eqno(2.2)
$$}

\medskip

As shown in [11] Prop. 2.5 we could assume only that the domain of the
operators $L_\ell$ contains a vector space $D$ which is a core for $G$
and that (2.2) holds for all $u,v\in D$.

For all $X\in{\cal B}(h)$ consider the sesquilinear form ${\cal L}(X)$
in $h$ with domain $D(G)\times D(G)$ given by
$$
\left\langle v,{\cal L}(X)u\right\rangle =
\left\langle v,XGu\right\rangle +
\left\langle Gv,Xu\right\rangle +
\sum_{\ell=1}^\infty\left\langle L_\ell v,XL_\ell u\right\rangle\eqno(2.3)
$$
Under the assumption {\bf A} it is possible to construct a q.d.s.
${\cal T}$ satisfying the equation
$$
\left\langle v,{\cal T}_t(X)u\right\rangle =
\left\langle v,Xu\right\rangle +
\int_0^t\left\langle v,
{\cal L}\left({\cal T}_s(X)\right)u\right\rangle ds\eqno(2.4)
$$
for all $u,v\in D(G)$  and all $X\in {\cal B}(h)$. As a first step
one proves the following

\proclaim Proposition {2.3}. Suppose that condition {\bf A}
holds and, for all $X\in {\cal B}(h)$, let $({\cal T}_t(X))_{t\ge 0}$ be
a strongly continuous family of elements of ${\cal B}(h)$ satisfying
(2.1). The following conditions are equivalent:
\item{\rm (i)} equation (2.4) holds for all $v,u\in D(G)$,
\item{\rm (ii)} for all $v,u\in D(G)$  we have
$$
\left\langle v,{\cal T}_t(X)u\right\rangle =
\left\langle P(t)v,XP(t)u\right\rangle +
\sum_{\ell=1}^\infty\int_0^t\left\langle L_\ell P(t-s)v,
{\cal T}_s(X)L_\ell P(t-s)u\right\rangle ds. \eqno(2.5)
$$

\noindent{\bf Proof.}\quad In order to show that condition i implies
condition ii we fix $t$ and compute the derivative
$$
{d\over ds}\left\langle P(t-s)v,{\cal T}_s(X)P(t-s)u\right\rangle =
\sum_{\ell=1}^\infty \left\langle L_\ell P(t-s)v,
{\cal T}_s(X) L_\ell P(t-s)u\right\rangle
$$
using equation (2.4). Clearly (2.5) follows integrating this identity
on the interval $[0,t]$. We prove now that condition ii implies
condition i. We recall first that $D(G^2)$ is a core for $G$ and,
for all $v,u\in D(G^2)$, we compute the derivative
$$\eqalign{
{d\over dt}\left\langle v,{\cal T}_t(X)u\right\rangle
&= \left\langle P(t)v,XP(t)Gu\right\rangle
+\left\langle P(t)Gv,XP(t)u\right\rangle\cr
& +\sum_{\ell=1}^\infty \left\langle L_\ell v,
{\cal T}_t(X)L_\ell u\right\rangle \cr
&+\int_0^t \sum_{\ell=1}^\infty\left\langle L_\ell P(t-s)Gv,
{\cal T}_s(X)L_\ell P(t-s)u\right\rangle ds\cr
&+\int_0^t \sum_{\ell=1}^\infty\left\langle L_\ell P(t-s)v,
{\cal T}_s(X)L_\ell P(t-s)Gu\right\rangle ds.\cr}
$$
The right-hand side of the above equation coincides with
$\left\langle v,{\cal L}({\cal T}_t(X))u\right\rangle$ by (2.4).
Therefore (2.4), for $v,u\in D(G^2)$, follows by integration on
$[0,t]$. Since $D(G^2)$ is a core for $G$ the proof is complete.
\quad\qed

\medskip

A solution of the equation (2.5) is obtained by the iterations
$$\eqalign{
\left\langle u, {\cal T}_t^{(0)}(X)u\right\rangle &=
\left\langle P(t)u,XP(t)u\right\rangle \cr
\left\langle u,{\cal T}_t^{(n+1)}(X)u\right\rangle &=
\left\langle P(t)u,XP(t)u\right\rangle \cr
&+\sum_{\ell=1}^\infty\int_0^t\left\langle L_\ell P(t-s)u,
{\cal T}_s^{(n)}(X)L_\ell P(t-s)u\right\rangle ds\cr} \eqno(2.6)
$$
(with $u\in D)$. Indeed, for all positive elements $X$ of
${\cal B}(h)$ and all $t\ge 0$, the sequence of operators
$\left({\cal T}_t^{(n)}(X)\right)_{n\ge 0}$ is non-decreasing.
Therefore it is strongly convergent and its limits for $X\in {\cal B}(h)$
and $t\ge 0$ define the {\sl minimal solution} ${\cal T}^{\rm (min)}$
of (2.5). We refer to [7], [12], [13] for more details.
The name ``minimal" is justified by the fact that,  given another
solution $({\cal T}_t)_{t\ge 0}$ of (2.4), one can easily prove that
$$
{\cal T}^{\rm (min)}_t(X)\le {\cal T}_t(X)\le \Vert X\Vert I\eqno(2.7)
$$
for all positive elements $X$ of ${\cal B}(h)$ and all $t\ge 0$
(see [7], [13]).
The minimal solution however, in spite of (2.2), is possibly
non-conservative  (see e.g. [12] Ex.~3.3 p.~174, [6] Ex.~3.6, 3.7 p. 97).

The infinitesimal generator ${\cal L}^{\rm (min)}$ will be given by
$$
{\cal L}^{\rm (min)}(X) = {\rm w}^*-\lim_{t\to 0^+}
t^{-1}\left({\cal T}^{\rm (min)}_t(X)-X\right)
$$
for all $X\in {\cal B}(h)$ such that the limit exists in the
ultraweak topology on ${\cal B}(h)$. Hence ${\cal L}$ is an extension
of ${\cal L}^{\rm (min)}$. Moreover it can be shown that the semigroup
${\cal T}^{\rm (min)}$ is conservative if and only if the identity
operator $I$ belongs to the domain of ${\cal L}^{\rm (min)}$
and ${\cal L}^{\rm (min)}(I)=0$ (see e.g. [9]).

\vskip 0.6truecm
\centerline{\bf 3. A representation of the resolvent of the minimal q.d.s.}

Let us consider the linear monotone maps
${\cal P}_\lambda:{\cal B}(h)\to {\cal B}(h)$ and
${\cal Q}_\lambda:{\cal B}(h)\to {\cal B}(h)$ defined by
$$\eqalignno{
\left\langle v,{\cal P}_\lambda(X)u\right\rangle
&=\int_0^{\infty}
\exp(-\lambda s)\left\langle P(s)v,X P(s)u\right\rangle\,ds &(3.1)\cr
\left\langle v,{\cal Q}_\lambda(X)u\right\rangle
&=\sum_{\ell=1}^\infty\int_0^{\infty}
\exp(-\lambda s)\left\langle L_\ell P(s)v,
X L_\ell P(s)u\right\rangle\,ds &(3.2)\cr}
$$
for all $\lambda>0$ and $X\in{\cal B}(h)$, $v,u\in D(G)$. It is easy to
check that both ${\cal P}_\lambda$ and ${\cal Q}_\lambda$ are completely
positive and normal contractions in ${\cal B}(h)$ (see e.g. [11] sect.2).

The resolvent of the minimal q.d.s.
$({\cal R}^{\rm (min)}_\lambda)_{\lambda>0}$ defined by
$$
\left\langle v,{\cal R}^{\rm (min)}_\lambda(X) u\right\rangle
= \int_0^\infty\exp(-\lambda s)
\left\langle v,{\cal T}^{\rm (min)}_s(X)u\right\rangle ds
$$
(with $X\in{\cal B}(h)$ and $v,u\in h$) can be represented as follows:

\proclaim Theorem {3.1}. For every $\lambda>0$ and $X\in {\cal B}(h)$
we have
$$
{\cal R}^{\rm (min)}_\lambda(X) = \sum_{k=0}^\infty
{\cal Q}^k_\lambda\left({\cal P}_\lambda(X)\right) \eqno(3.3)
$$
the series being convergent for the strong operator topology.

\noindent{\bf Proof.}\quad Consider the sequence
$({\cal R}^{(n)}_\lambda)_{n\ge 0}$ of linear monotone maps
${\cal R}^{(n)}_\lambda:{\cal B}(h)$ $\to {\cal B}(h)$ given by
$$
\left\langle v,{\cal R}^{(n)}_\lambda(X)u\right\rangle
=\int_0^\infty\exp(-\lambda s)
\left\langle v,{\cal T}^{(n)}_s(X)u\right\rangle ds
$$
where the operators ${\cal T}^{(n)}_s$ are defined by (2.6).
Clearly (2.1) guarantees that ${\cal R}^{(n)}_\lambda$ is well defined.
Moreover, for all positive elements $X$ of ${\cal B}(h)$, the
sequence $\left({\cal R}^{(n)}_\lambda(X)\right)_{n\ge 0}$ is
non-decreasing. Therefore, by the definition of minimal q.d.s., for
all $u\in h$ we have
$$
\left\langle u,{\cal R}^{\rm (min)}_\lambda(X)u\right\rangle
=\sup_{n\ge 0} \left\langle u,{\cal R}^{(n)}_\lambda(X)u\right\rangle
=\int_0^\infty\exp(-\lambda s)
\left\langle u,{\cal T}^{\rm (min)}_s(X)u\right\rangle ds.
$$
The second equation (2.6) yields
$$\eqalign{
\left\langle u,{\cal R}^{(n+1)}_\lambda(X)u\right\rangle
&=\int_0^\infty \hbox{e}^{-\lambda t}
\left\langle P(t)v,XP(t)u\right\rangle dt\cr
& \ + \sum_{\ell=1}^\infty\int_0^\infty \hbox{e}^{-\lambda t} dt
\int_0^t \left\langle L_\ell P(t-s)u,{\cal T}^{(n)}_s(X)
L_\ell P(t-s)u\right\rangle ds\cr}
$$
for all $u,v\in D(G)$. By the change of variables $(r,s)=(t-s,s)$ in
the above double integral we have
$$\eqalign{
\left\langle u,{\cal R}^{(n+1)}_\lambda(X)u\right\rangle
&=\left\langle u,{\cal P}_\lambda(X) u\right\rangle\cr
& \ + \sum_{\ell=1}^\infty\int_0^\infty \hbox{e}^{-\lambda r} dr
\int_0^\infty \hbox{e}^{-\lambda s}
\left\langle L_\ell P(r)u,{\cal T}^{(n)}_s(X)
L_\ell P(r)u\right\rangle ds.\cr}
$$
Thus we obtain the recursion formula
$$
{\cal R}^{(n+1)}_\lambda(X) = {\cal P}_\lambda(X)
+ {\cal Q}_\lambda({\cal R}^{(n)}_\lambda(X)).
$$
Iterating $n$ times this equation we have
$$
{\cal R}^{(n+1)}_\lambda(X) =
\sum_{k=0}^{n+1}{\cal Q}^k_\lambda ({\cal P}_\lambda(X))
$$
and (3.3) follows letting $n$ tend to $+\infty$. Clearly (3.3) also holds
for an arbitrary element of ${\cal B}(h)$ since each bounded operator can
be written as a linear combination of positive self-adjoint operators.
\quad\qed

\medskip

The following proposition gives another useful relation between
${\cal P}_\lambda$, ${\cal Q}_\lambda$ and
${\cal R}^{\rm (min)}_\lambda$.

\proclaim Proposition {3.2}. Suppose that the condition {\bf A} holds
and fix $\lambda>0$. For all $n\ge 1$ we have
$$
\sum_{k=0}^n {\cal Q}^{k}_\lambda({\cal P}_\lambda(I))
+\lambda^{-1}{\cal Q}^{n+1}_\lambda(I) = \lambda^{-1}I.\eqno(3.4)
$$

\noindent{\bf Proof.}\quad For all $u\in D(G)$ a straightforward
computation yields
$$\eqalign{
\sum_{\ell=1}^\infty \int_0^\infty
\hbox{e}^{-\lambda t}\left\Vert L_\ell P(t)u\right\Vert^2 dt
&= -2\re\int_0^\infty \hbox{e}^{-\lambda t}
\left\langle  P(t)u,GP(t)u\right\rangle dt\cr
&= -\int_0^\infty \hbox{e}^{-\lambda t}
{d\over dt} \left\Vert P(t)u\right\Vert^2 dt \cr
& = \left\Vert u\right\Vert^2 -\lambda
\int_0^\infty \hbox{e}^{-\lambda t}\left\Vert P(t)u\right\Vert^2 dt
\cr} \eqno(3.5)
$$
Therefore we have
$$
{\cal P}_\lambda(I) + \lambda^{-1}{\cal Q}_\lambda(I)
=\lambda^{-1}(I).
$$
This proves (3.4) for $n=0$. Suppose it has been established for an
integer $n$. Applying the map ${\cal Q}_\lambda$ to both sides of (3.4)
yields
$$
\sum_{k=1}^{n+1}{\cal Q}^k_\lambda({\cal P}_\lambda(I))
+\lambda^{-1}{\cal Q}^{n+2}_\lambda(I) = \lambda^{-1}{\cal Q}_\lambda(I)
=\lambda^{-1}I - {\cal P}_\lambda(I).
$$
This proves (3.4) for the integer $n+1$ and completes the proof.\quad\qed

\medskip

The representation formula for the resolvent of the minimal q.d.s.
of Theorem 3.1 allows us to give a quick proof of the following necessary
and sufficient condition for conservativity obtained by the first author
in [7] (see also [11] Prop. 2.7).

\proclaim Proposition {3.3}. Suppose that the condition {\bf A} holds
and fix $\lambda>0$. Then the sequence of positive operators
$({\cal Q}^n_\lambda(I))_{n\ge 0}$ is non-increasing.
Moreover the following conditions are equivalent:
\item{\rm (i)} the q.d.s. ${\cal T}^{\rm (min)}$ is conservative,
\item{\rm (ii)} s$-\lim_{n\to\infty}{\cal Q}_\lambda^n(I) =0$.

\noindent{\bf Proof.}\quad The sequence of positive operators
$({\cal Q}^n_\lambda(I))_{n\ge 0}$ is non-increasing because
(3.4) yields
$$
{\cal Q}^{n}_\lambda(I) - {\cal Q}^{n+1}_\lambda(I)
=\lambda {\cal Q}^{n}_\lambda({\cal P}_\lambda(I)).
$$
Therefore it is strongly convergent to a positive operator $Y$.
Letting $n$ tend to $+\infty$ in (3.4), we have
$$
{\cal R}^{\rm (min)}_\lambda(I) +\lambda^{-1}Y
= \lambda^{-1} I.
$$
Now condition i can be clearly stated as:
${\cal R}^{\rm (min)}_\lambda(I)$ coincides with $\lambda^{-1}I$.
Therefore the desired equivalence follows.\quad\qed

\vskip 0.6truecm
\centerline{\bf 4. Sufficient conditions for conservativity}

The minimal q.d.s. is conservative whenever, for a fixed $\lambda >0$,
the series
$$
\sum_{k=1}^\infty\left\langle u,{\cal Q}_\lambda^{k}(I)u\right\rangle
\eqno(4.1)
$$
is convergent for all $u$ in a dense subspace of $h$.
In fact, in this case, condition (ii) of Proposition 3.3 holds
because the sequence of positive operators
$({\cal Q}^{k}_\lambda(I))_{k\ge 1}$ is non-increasing.

In this section we shall give easily verifiable conditions on the
operators $G$, $L_\ell$ that guarantee convergence of (4.1).
As a first step we shall prove an easy estimate.
Let $R(n;G)$ be the resolvent operator $(nI-G)^{-1}$.
The operator in $h$ with domain $D(G)$
$$
\sum_{\ell=1}^\infty
\left(nL_\ell R(n;G)\right)^*\left(nL_\ell R(n;G)\right)
$$
has a unique bounded extension by virtue of identity (2.2) and
well-known properties of resolvent operators. We shall denote this
bounded extension by $F_n$. Notice that $F_n$ is a positive
self-adjoint operator.

\proclaim Proposition {4.1}. For every $u\in h$ we have
$$
\sum_{k=1}^\infty\left\langle u,{\cal Q}_\lambda^{k}(I)u\right\rangle
\le\liminf_{n\to\infty}
\left\langle u,{\cal R}^{\rm (min)}_\lambda(F_n)u\right\rangle.
$$

\noindent{\bf Proof.}\quad  For $u\in D(G)$, $n\ge 1$, we have
$$\eqalign{
\left\langle u,{\cal P}_\lambda(F_n)u\right\rangle
&=\sum_{\ell=1}^\infty\int_0^\infty\hbox{e}^{-\lambda t}
\left\Vert nL_\ell R(n;G)P(t)u\right\Vert^2dt\cr
&=\sum_{\ell=1}^\infty\int_0^\infty\hbox{e}^{-\lambda t}
\left\Vert L_\ell P(t)(nR(n;G)u)\right\Vert^2 dt\cr
&=\left\langle nR(n;G)u,n{\cal Q}_\lambda(I)R(n;G)u\right\rangle.\cr}
$$
Therefore the bounded operators ${\cal P}_\lambda(F_n)$ and
$n^2R(n;G^*){\cal Q}_\lambda(I)R(n;G)$ coincide. Moreover the sequence
of operators $({\cal P}_\lambda(F_n))_{n\ge 1}$ is uniformly bounded
and converges strongly to ${\cal Q}_\lambda(I)$ by well-known properties
of resolvent operators.

The maps ${\cal Q}_\lambda^k$ are normal (cf. Definition 2.1 v).
We have then, for $u\in h$,
$$
\sum_{k=1}^\infty\left\langle u,{\cal Q}_\lambda^{k}(I)u\right\rangle
\le\liminf_{n\to\infty}\sum_{k=0}^\infty
\left\langle u,{\cal Q}_\lambda^k({\cal P}_\lambda(F_n))u\right\rangle
= \liminf_{n\to\infty}\left\langle u,
{\cal R}^{\rm (min)}_\lambda(F_n)u\right\rangle,
$$
by Fatou's lemma and Theorem 3.1.\quad\qed

\medskip

In order to estimate ${\cal R}^{\rm (min)}_\lambda(C)$
for self-adjoint operators $C$  we introduce now our key assumption

\smallskip
{\bf C} - {\sl A positive self-adjoint operator $C$ satisfies
Condition {\bf C} if:

--  the domain of its positive square root $C^{1/2}$ contains the domain
$D(G)$ of $G$ and $D(G)$ is a core for $C^{1/2}$,

-- the linear manifolds $L_\ell(D(G^2))$,
$\ell\ge 1$, are  contained in the domain of $C^{1/2}$,

-- there exists a positive constant $b$ such that
$$
2\re \left\langle C^{1/2}u,C^{1/2}Gu\right\rangle
+\sum_{\ell=1}^\infty\left\langle C^{1/2}L_\ell u,
C^{1/2}L_\ell u\right\rangle
\le b \left\Vert C^{1/2}u\right\Vert^2 \eqno(4.2)
$$
for all $u\in D(G^2)$.}
\smallskip

\noindent{\bf Remark.} Condition {\bf C} implies that, for each
$u\in D(G^2)$, the function $t\to\left\Vert C^{1/2}P(t)u\right\Vert^2$
is differentiable  and
$$
{d\over dt} \left\Vert C^{1/2}P(t)u\right\Vert^2
= 2\re\left\langle C^{1/2}P(t)u,C^{1/2}GP(t)u\right\rangle.
$$
Indeed, for each $u\in D(G)$ and each $\lambda>0$, let
$v=\lambda^{-1}R(\lambda;G)u$. The inequality (4.2) yields
$$\eqalign{
\left\Vert C^{1/2}u\right\Vert^2
& = \left\Vert C^{1/2}v\right\Vert^2
- 2\lambda^{-1}\re\left\langle C^{1/2}v,C^{1/2}Gv\right\rangle
+ \lambda^{-2}\left\Vert C^{1/2}Gv\right\Vert^2\cr
& \ge \left(1-\lambda^{-1}b\right)\left\Vert C^{1/2}v\right\Vert^2
= \left( 1-\lambda^{-1}b\right)
\left\Vert C^{1/2}\lambda R(\lambda;G)u\right\Vert^2. \cr}
$$
The above inequality also holds for $u\in D(C^{1/2})$ since $D(G)$
is a core for $C^{1/2}$. It follows that $G$ is the
infinitesimal generator of a strongly continuous semigroup on
the Hilbert space $D(C^{1/2})$ (endowed with the graph norm).
This is obtained by restricting the operators $P(t)$ to $D(C^{1/2})$.
Therefore the claimed differentiation formula follows.

\smallskip
Under assumption {\bf C} we can prove a useful estimate of
${\cal R}^{\rm (min)}_\lambda(C_\varepsilon)$ where
$(C_\varepsilon)_{\varepsilon>0}$ is the family of bounded
regularization $C_\varepsilon=C(I+\varepsilon C)^{-1}$.

\proclaim Proposition {4.2}. Suppose that conditions {\bf A} and {\bf C}
hold. Then, for all $\lambda>b$ and all $u\in D(G^2)$, we have
$$
(\lambda-b)\sup_{\varepsilon > 0}\left\langle u,
{\cal R}^{\rm (min)}_\lambda(C_\varepsilon)u\right\rangle
\le \left\Vert C^{1/2}u\right\Vert^2.\eqno(4.3)
$$

\noindent{\bf Proof.}\quad
Let $\left({\cal R}^{(n)}_\lambda\right)_{n>0}$ be the sequence of
monotone linear maps considered in the proof of Theorem 3.1. Clearly it
suffices to show that, for all $n\ge 0$, $\lambda > b$ and $u\in D(G^2)$,
the operator ${\cal R}^{(n)}_\lambda(C_\varepsilon)$ satisfies
$$
(\lambda-b)\sup_{\varepsilon>0}\left\langle u,
{\cal R}^{(n)}_\lambda(C_\varepsilon)u\right\rangle
\le \left\Vert C^{1/2}u\right\Vert^2.\eqno(4.4)
$$
The above inequality holds for $n=0$. In fact, integrating by parts,
for all $u\in D(G^2)$, we have
$$\eqalign{
\lambda\left\langle u,
{\cal R}^{(0)}_\lambda(C_\varepsilon)u\right\rangle
& = \lambda\int_0^\infty\hbox{e}^{-\lambda t}
\left\langle P(t)u,C_\varepsilon P(t)u\right\rangle dt\cr
& \le \lambda\int_0^\infty\hbox{e}^{-\lambda t}
\left\Vert C^{1/2} P(t)u\right\Vert^2 dt\cr
& =\left\Vert C^{1/2}u\right\Vert^2
+ 2\re \int_0^\infty\hbox{e}^{-\lambda t}
\left\langle C^{1/2}P(t)u, C^{1/2}GP(t)u\right\rangle dt.
\cr}
$$
The inequality (4.2) yields
$$\eqalign{
\lambda\left\langle u,
{\cal R}^{(0)}_\lambda(C_\varepsilon)u\right\rangle
& \le \left\Vert C^{1/2}u\right\Vert^2
+ b \int_0^\infty \hbox{e}^{-\lambda t}
\left\Vert C^{1/2}P(t)u\right\Vert^2 dt\cr
& = \left\Vert C^{1/2}u\right\Vert^2
+ b \sup_{\varepsilon>0}
\left\langle u,{\cal R}^{(0)}_\lambda(C_\varepsilon)u\right\rangle.\cr}
$$
This clearly implies (4.4) for $n=0$. Suppose that (4.4) has been established
for an  integer $n$; then, from the second equation (2.6) and the definition
of ${\cal R}^{(n)}_\lambda$, we have
$$\eqalign{
\left\langle u,
{\cal R}^{(n+1)}_\lambda(C_\varepsilon)u\right\rangle
&=\left\langle u,{\cal P}_\lambda(C_\varepsilon)u\right\rangle\cr
&+\sum_{l=1}^\infty \int_0^\infty\hbox{e}^{-\lambda t}
\left\langle L_\ell P(t)u,
{\cal R}^{(n)}_\lambda(C_\varepsilon)L_\ell P(t)u\right\rangle dt \cr
&\le \left\langle u,{\cal P}_\lambda(C_\varepsilon)u\right\rangle
+ {1\over \lambda-b}
\sum_{\ell=1}^\infty\int_0^\infty \hbox{e}^{-\lambda t}
\left\Vert C^{1/2}L_\ell P(t)u\right\Vert^2 dt
\cr}
$$
Inequality (4.2) and integration by parts yield
$$\eqalign{
\sum_{\ell=1}^\infty\int_0^\infty\kern-1.6truept\hbox{e}^{-\lambda t}
\left\Vert C^{1/2}L_\ell P(t)u\right\Vert^2 dt
& \le \int_0^\infty \kern-1.6truept\hbox{e}^{-\lambda t}
\left(-{d\over dt}\left\Vert C^{1/2} P(t)u\right\Vert^2\right)dt\cr
&+ b\int_0^\infty \kern-1.6truept\hbox{e}^{-\lambda t}
\left\Vert C^{1/2} P(t)u\right\Vert^2 dt\cr
&=\left\Vert C^{1/2}u\right\Vert^2
-(\lambda-b)\int_0^\infty \kern-1.6truept\hbox{e}^{-\lambda t}
\left\Vert C^{1/2} P(t)u\right\Vert^2 dt\cr
&\le\left\Vert C^{1/2}u\right\Vert^2
-(\lambda-b)\left\langle u,{\cal P}_\lambda(C_\varepsilon)u\right\rangle.\cr}
$$
Therefore (4.4) for $n+1$ follows. The proof is complete.\quad\qed

\medskip

We can now prove the main results of this paper.

\proclaim Theorem {4.3}. Suppose that condition {\bf A} holds and
there exists an operator $C$ satisfying condition {\bf C} such that
$$
\left\langle u,F_n u\right\rangle \le \left\langle u,Cu\right\rangle
$$
for all $u\in D(C)$, $n\ge 1$. Then the minimal q.d.s. is conservative.

\noindent{\bf Proof.}\quad Let $\lambda>b$ fixed. Under the
present hypotheses, for $\varepsilon >0$, the bounded operators
$(F_n)_\varepsilon$ and $C_\varepsilon$  satisfy the inequality
$(F_n)_\varepsilon \le C_\varepsilon$ (see, e.g. [22] Chap. 8, Ex. 51, p.317).
Applying Proposition 4.2, we obtain the estimate
$$\eqalign{
\sum_{k=1}^\infty\left\langle u,{\cal Q}_\lambda^{k}(I)u\right\rangle
&\le \liminf_{n\to\infty}\sup_{\varepsilon>0}
\left\langle u,
{\cal R}^{\rm (min)}_\lambda((F_n)_\varepsilon)u\right\rangle\cr
&\le\sup_{\varepsilon>0}
\left\langle u,{\cal R}^{\rm (min)}_\lambda(C_\varepsilon)u\right\rangle
\le (\lambda-b)^{-1}\left\Vert C^{1/2}u\right\Vert^2 < +\infty\cr}
$$
Therefore the minimal q.d.s. is conservative since
condition (b) of Proposition 3.3 is fulfilled.\quad\qed

\medskip

Notice that, in the above theorem, we did not assume that the
quadratic form
$$
u\to -2\re\left\langle u,Gu\right\rangle
$$
with domain $D(G)$ is closable (see Example 5.3).

\proclaim Theorem {4.4}. Suppose that assumptions {\bf A}, {\bf C} hold
for some positive self-adjoint operator $C$ and there
exists a positive self-adjoint operator $\Phi$ in $h$ such that:
\item{\rm (a)} the domain of the positive square root $\Phi^{1/2}$
contains the domain of $G$ and, for every $u\in D(G)$, we have
$$
- 2\re\left\langle u,Gu\right\rangle
=\sum_{\ell=1}^\infty \left\langle L_\ell u,L_\ell u\right\rangle
=\left\langle \Phi^{1/2}u,\Phi^{1/2}u\right\rangle,
$$
\item{\rm (b)} the domain of $C$ is contained in the domain
of $\Phi$ and, for every $u\in D(C)$, we have
$$
\left\langle \Phi^{1/2}u,\Phi^{1/2}u\right\rangle
\le \left\langle C^{1/2}u,C^{1/2}u\right\rangle.
$$
Then the minimal q.d.s. is conservative.

\noindent{\bf Proof.}\quad Let $\lambda>b$ and $u\in D(G^2)$ fixed.
For $\varepsilon >0$, the bounded operators $\Phi_\varepsilon$
and $C_\varepsilon$ satisfy the inequality
$\Phi_\varepsilon \le C_\varepsilon$
(see, e.g. [22] Chap. 8, Ex. 51, p.317).  Moreover, for $u\in D(G)$,
we have
$$\eqalign{
\sup_{\varepsilon>0}
\left\langle u,{\cal P}_\lambda(\Phi_\varepsilon)u\right\rangle
&=\int_0^\infty\hbox{e}^{-\lambda t}
\left\Vert \Phi^{1/2}P(t)u\right\Vert^2dt\cr
&=\sum_{\ell=1}^\infty\int_0^\infty\hbox{e}^{-\lambda t}
\left\Vert L_\ell P(t)u\right\Vert^2 dt
=\left\langle u,{\cal Q}_\lambda(I)u\right\rangle.\cr}
$$
This implies that the non-decreasing family of operators
$({\cal P}_\lambda(\Phi_\varepsilon))_{\varepsilon>0}$
is uniformly  bounded and, since $D(G)$ is dense in $h$,
it follows that it converges strongly to ${\cal Q}_\lambda(I)$
as $\varepsilon$ goes to $0$. The maps ${\cal Q}_\lambda^k$ being
normal we have
$$
\sum_{k=0}^\infty\left\langle u,{\cal Q}_\lambda^{k+1}(I)u\right\rangle
=\sup_{\varepsilon >0}\sum_{k=0}^\infty
\left\langle u,{\cal Q}^k({\cal P}_\lambda(\Phi_\varepsilon))u\right\rangle
=\sup_{\varepsilon >0}\left\langle u,
{\cal R}^{\rm (min)}_\lambda(\Phi_\varepsilon)u\right\rangle
$$
by Theorem 3.1. Applying Proposition 4.2 we obtain the estimate
$$\eqalign{
\sum_{k=1}^\infty\left\langle u,{\cal Q}_\lambda^{k}(I)u\right\rangle
&=\sup_{\varepsilon >0}\left\langle u,
{\cal R}^{\rm (min)}_\lambda(\Phi_\varepsilon)u\right\rangle\cr
&\le \sup_{\varepsilon >0}\left\langle u,
{\cal R}^{\rm (min)}_\lambda(C_\varepsilon)u\right\rangle
\le (\lambda-b)^{-1}\left\Vert C^{1/2}u\right\Vert^2\cr}
$$
Therefore the minimal q.d.s. is conservative because
condition (b) of Proposition 3.3 holds.\quad\qed

\medskip

The following corollary gives a simpler and easily verifiable condition
under stronger assumptions on the domain of the operator $C$.

\proclaim Corollary {4.5}. Suppose that assumption {\bf A}
holds and there exist a self-adjoint operator $C$ and a core $D$
for $G$ with the following properties:
\item{(a)} the domain of $G$ coincides with the domain of $C$ and for
all $u\in D(G)$ there exists a sequence $(u_n)_{n\ge 0}$ of elements
of $D$ such that both $(Gu_n)_{n\ge 0}$ and $(Cu_n)_{n\ge 0}$ converge
strongly,
\item{(b)} there exists a positive self-adjoint operator $\Phi$ such that
the domain of $\Phi$ contains the domain of $D$ and for all $u\in D$ and
$n\ge 1$ we have the inequality
$$
-2\re\left\langle u,Gu\right\rangle
= \left\langle u,\Phi u\right\rangle
\le \left\langle u,Cu\right\rangle,
$$
\item{(c)} for all $\ell\ge 1$, $L_\ell(D)\subseteq D(C)$,
\item{(d)} there exists a constant $b$ such that, for all $u\in D$,
the following inequality holds
$$
2\re \left\langle Cu,Gu\right\rangle
+\sum_{\ell=1}^\infty \left\langle L_\ell u,CL_\ell u\right\rangle
\le b\left\langle u, Cu\right\rangle.\eqno(4.5)
$$
Then the minimal q.d.s. is conservative.

\noindent{\bf Proof.}\quad The inequality of condition (b) obviously
holds also for $u\in D(C)$ because of condition (a) and self-adjointness
of $\Phi$. Therefore, in order to prove the corollary, it suffices
to show that, under the above hypotheses, the operator $C$
satisfies assumption {\bf C} and apply Theorem 4.4.

Let $(u_n)_{n\ge 0}$ be a sequence of elements of $D$ such that
$$
\lim_{n\to\infty}Cu_n = Cu,\qquad\qquad
\lim_{n\to\infty}Gu_n = Gu.
$$
Condition (d) implies that $(C^{1/2}L_\ell u_n)_{n\ge 1}$ is a
Cauchy sequence for $\ell\ge 1$. Therefore it is convergent
and it is easy to deduce that (4.2) holds for $u\in D(G)$. \quad\qed

\medskip

\noindent{\bf Remark.}\quad Another simple sufficient condition
for conservativity can be easily obtained by substituting (b) in
the above corollary with the following hypothesis:
\item{(b')} for all $u\in D$ and $n\ge 1$ we have the inequality
$$
-2\re\left\langle nR(n;G)u,nGR(n;G)u\right\rangle
= \left\langle u,F_n u\right\rangle
\le \left\langle u,Cu\right\rangle.
$$
The proof can be easily done by applying Theorem 4.3.

\vskip 0.6truecm
\centerline{\bf 5. Applications and examples}

In this section we apply our results to study conservativity of
three minimal q.d.s.: one arising from a physical model and another
>from extension problems of classical Markovian semigroups to
non-abelian algebras. We consider semigroups of diffusion type since
the minimal q.d.s. of jumps and drift type leaves the abelian algebra
of multiplication operators invariant. Therefore the conservativity
problem for the minimal q.d.s. can be reduced to a problem in classical
probability.

\medskip
\noindent{\bf 5.1 \quad Q.d.s. in a model for heavy ion collision.}
\par\noindent
As a first example we apply the conservativity condition of Corollary 4.5
to the minimal q.d.s. proposed by Alicki (see, for example, [3], [4])
to describe phenomenologically a quantum system with dissipative heavy
ion collisions. This problem can not be solved by applying the tools
developed in [9], [11].

Let $h=L^2(\IR^3;\IC)$ and let $m\in ]0,+\infty[$, $\alpha\in\IR$.
We denote by $\partial_\ell$ ($\ell=1,2,3$) differentiation with respect
to the $\ell$-th coordinate. Let $V:\IR^3\to\IR$, $W:\IR^3\to\IR$ be two
functions with the following  properties:
\item{1.} $V$ can be written as the sum of a bounded function and a
square integrable function, $V$ is differentiable and the partial
derivatives $\partial_\ell V$ are bounded,
\item{2.} $W$ is bounded and
$$
\sup_{x\in\IR^3}|W(x)|^2 < \left(m\alpha^2\right)^{-1},
$$
$W$ is twice differentiable and the following functions
are bounded continuous
$$
x\to x_\ell W(x),\qquad
x\to x_\ell\partial_\ell W(x),\qquad
x\to \partial_\ell^2 W(x)\quad \ell=1,2,3. \eqno(5.1)
$$

Consider the operators $H_0$, $V$, $L_\ell$, $G$ with domain $H^2(\IR;\IC)$
$$\eqalign{
H_0u = -{1\over 2m}\Delta u,\qquad   (Vu)(x) & = V(x)u(x),\cr
L_\ell u = W(x)(x_\ell + \alpha\partial_\ell)u,\qquad
Gu & = -i(H_0 +V)u -{1\over 2}\sum_{\ell=1}^3 L^*_\ell L_\ell u\cr}
$$
and let $L_\ell=0$ for $\ell\ge 4$.
The arguments of [19] Ch. V, Sect. 3 show that the operator $G$ is
a relatively bounded perturbation of $H_0$ with relative bound
smaller than $1$ and the linear manifold $D$ of infinitely
differentiable functions  with compact support is a core for $G$.
The operator $G$ is the infinitesimal generator of a strongly continuous
contraction semigroup in $h$ by [19] Th. 2.7 p. 499 and following remarks.
Thus the basic assumption {\bf A} holds because it suffices to check
identity (2.2) for $v,u\in D$  and, in this case, (2.2) is trivial.

We show that the minimal q.d.s. is conservative applying Corollary 4.5.
The most natural choice of the operator $C$ is the following
$$
D(C)=H^2(\IR^3;\IC),\qquad\quad Cu = c \left(-\Delta + 1\right)u
$$
where $c$ is a suitable constant to be determined. In fact
hypothesis (a) obviously holds because $G$ and $C$ are relatively
bounded one respect to the other and $D$ is a core for both.
By virtue of von Neumann's theorem (see [19] Th. 3.24 p. 275)
(b) is satisfied. Hypothesis (c) is trivially fulfilled.
In order to check (d) notice first that it suffices to check (4.5) for
all $u\in D$ i.e. to estimate the quadratic form associated with
the formal operator
$$
CG + G^*C + \sum_{\ell=1}^3 L_\ell^* C L_\ell = i[ V,C ]
+{1\over 2}\sum_{\ell=1}^3
\left(L_\ell^*[C,L_\ell] + [L_\ell^*,C]L_\ell\right).
$$
This turns out to be a second order differential operator with bounded
coefficients. Hence, for $u\in D$, we have
$$
2\re \left\langle Cu,Gu\right\rangle
+ \sum_{\ell=1}^\infty \left\langle L_\ell u, CL_\ell u\right\rangle
\le b \left\langle u, Cu\right\rangle
$$
where $b$ depends only on the supremum of the partial derivatives
$\partial_\ell V$, $\ell=1,2,3$, and of the functions (5.1). Therefore
hypothesis (d) of Corollary (4.5) also holds and the minimal q.d.s.
is conservative.

>From the above discussion it is clear that our result can be applied
to a large class of Lindblad type perturbations of pure hamiltonian
evolutions arising in physical models (see, for instance, [4]).

\medskip
\noindent{\bf 5.2 Extension of classical Brownian motions with partial
reflection.}
\par\noindent
Let $h=L^2((0,+\infty);\IC)$, let $\alpha \in ]0,+\infty[$
and let $g$ be a function in $h$. Define the parameter
$\theta = \left\Vert g\right\Vert^2/(2\alpha)$. Consider the
operators  $G$ and $L_\ell$
$$\eqalign{
D(G)=\left\{ u\in H^2((0,+\infty);\IC)\mid u'(0) = \theta u(0)\right\}
\qquad Gu & = {1\over 2} u'',\cr
D(L_1)= H^1((0,+\infty);\IC)\hskip 5truecm L_1u & = u',\cr
D(L_2)= H^1((0,+\infty);\IC)
\hskip 5truecm L_2u & = {u(0)\over\sqrt{2\alpha}}g,\cr}
$$
and let $L_\ell =0$ for all $\ell\ge 3$. In [5] (Prop. 4.3 and Th. 2.4)
it has been shown that:
\item{1.} Our basic assumption {\bf A} is satisfied.
\item{2.} The operator $G$ is negative and self-adjoint.
\item{3.} The restriction of the map ${\cal L}$ defined by (2.3) to
multiplication operators by a regular bounded real function $f$ on
$[0,+\infty[$, coincides with the infinitesimal generator $A$ of a brownian
motion on $[0,+\infty[$ with partial reflection at the boundary point
$\{0\}$ and  partial reentrance in $]0,+\infty[$ with reentrance density
$x\to |g(x)|^2$.
$$\eqalign{
D(A) &=\left\{ f\in {\cal C}^2_b([0,+\infty[;\IR)\,\Big|\,
\alpha f'(0) + \int_0^\infty \left(f(x)-f(0)\right) |g(x)|^2 dx = 0
\right\}\cr
&\qquad\qquad\qquad\quad (Af)(x) = {1\over 2}f''(x)\cr}
$$
where ${\cal C}^2_b([0,+\infty[;\IR)$ denotes the vector space of
real functions on $[0,+\infty[$ bounded with bounded continuous
derivatives up to the second order.
\item{4.} The minimal q.d.s. is an extension to ${\cal B}(h)$ of the
Markovian semigroup of the classical stochastic process if and only if
it is conservative.

\smallskip
Here we apply the main result of this paper to show that the minimal
q.d.s. constructed from the above operators $G$ and  $L_\ell$ is
conservative whenever $g$ belongs to $H^1((0,+\infty);\IC)$.
We prove first the following

\proclaim Lemma {5.1}. For all $u\in H^1((0,+\infty);\IC)$
and $\lambda>0$ we have
$$
\lim_{\lambda\to +\infty} (\lambda R(\lambda;G)u)(0) = u(0),\qquad
s-\lim_{\lambda\to +\infty} \lambda L_1 R(\lambda;G)u = L_1u.
$$

\noindent{\bf Proof.}\quad  An elementary computation yields
$$\eqalign{
\left(\lambda R(\lambda;G)u\right)(x) &=
\sqrt{\lambda\over 2}
\int_0^{\infty}\exp\left(-\sqrt{2\lambda}|x-s|\right)u(s)ds\cr
&+\sqrt{\lambda\over 2} {\sqrt{2\lambda}-\theta\over\sqrt{2\lambda}+\theta}
\int_0^{\infty}\exp\left(-\sqrt{2\lambda}(x+s)\right)u(s)ds, \cr
\left(\lambda L_1R(\lambda;G)u\right)(x) &=
-\lambda\int_0^{\infty}\hbox{\rm sgn}(x-s)
\exp\left(-\sqrt{2\lambda}|x-s|\right)u(s)ds\cr
&-\lambda{\sqrt{2\lambda}-\theta\over\sqrt{2\lambda}+\theta}
\int_0^{\infty}\exp\left(-\sqrt{2\lambda}(x+s)\right)u(s)ds \cr}
$$
where $\hbox{\rm sgn}(x-s)=1$ if $x\ge s$ and $\hbox{\rm sgn}(x-s)=-1$
if $x<s$. Therefore the first limit is easily computed.
Integrating by parts both the above integrals we have
$$\eqalign{
\left(\lambda L_1R(\lambda;G)u\right)(x)
&=\left(\lambda R(\lambda;G)L_1u\right)(x)
+ {\theta\sqrt{2\lambda}\over\sqrt{2\lambda}+\theta }
\exp\left(-\sqrt{2\lambda} x\right) u(0)\cr
& - \sqrt{2\lambda}{\sqrt{2\lambda}-\theta\over\sqrt{2\lambda}+\theta}
\int_0^{\infty}\exp\left(-\sqrt{2\lambda}(x+s)\right)u'(s)ds\cr }
$$
The first term converges to the desired limit, for the strong topology
on $h$, by a well-known property of the resolvent operators and the
second clearly vanishes as $\lambda$ goes to $+\infty$.

Disregarding the factor
$(\sqrt{2\lambda}-\theta)/(\sqrt{2\lambda}+\theta)$
goes to $1$ as $\lambda$ goes to $+\infty $ and using the
Schwarz inequality we can estimate the third term by
$$\eqalign{
& 2\lambda \int_0^\infty  \exp\left(-2\sqrt{2\lambda}\,x\right)dx\cdot
\left|\int_0^\infty \exp\left(-\sqrt{2\lambda}s\right)
u'(s)ds\right|^2 \cr
& = {\sqrt{2\lambda}\over 2 }
\left|\int_0^\infty \exp\left(-\sqrt{2\lambda}s/2\right)\cdot
\left( \exp\left(-\sqrt{2\lambda}s/2\right)u'(s)\right)ds\right|^2\cr
 & \le {\sqrt{2\lambda}\over 2}
\int_0^\infty \exp\left(-\sqrt{2\lambda}s\right) ds
\cdot \int_0^\infty \exp\left(-\sqrt{2\lambda}s\right) |u'(s)|^2 ds \cr
& = {1\over 2}
\int_0^\infty \exp\left(-\sqrt{2\lambda}s\right) |u'(s)|^2 ds\cr}
$$
The right-hand side vanishes as $\lambda$ goes to $+\infty$ by
Lebesgue's theorem.

This completes the proof. \quad\qed

\proclaim Lemma {5.2}. Let $C$ be the positive self-adjoint operator $-2G$.
The vector space $H^1((0,+\infty);\IC)$ is contained in the domain of
$C^{1/2}$ and
$$
\left\Vert C^{1/2}u\right\Vert^2
=\left\Vert L_1u\right\Vert^2 +\left\Vert L_2u\right\Vert^2 \eqno(5.2)
$$
for all $u\in H^1((0,+\infty);\IC)$.

\noindent{\bf Proof.}\quad Let $u\in H^1((0,+\infty);\IC)$ and let
$u_\lambda = \lambda R(\lambda;G)u$. The vector $u_\lambda$ belongs
to the domains of $C$ and $G$. Moreover, for all $\lambda,\mu>0$,
(2.2) yields
$$
\left\Vert C^{1/2}(u_\lambda - u_\mu)\right\Vert^2
=\left\Vert L_1(u_\lambda - u_\mu)\right\Vert^2
+\left\Vert L_2(u_\lambda - u_\mu)\right\Vert^2.
$$
Therefore the family of vectors $(C^{1/2} u_\lambda)_{\lambda >0}$
is Cauchy by Lemma 5.1.  Thus $u$ belongs to the domain of $C^{1/2}$.
Moreover (5.2) holds for all vectors $u_\lambda$ because it is
equivalent to (2.2) for vectors belonging to $D(G)$. Letting $\lambda$
go to $+\infty$ we see that (5.2) holds also for the vector $u$.
\quad\qed

\medskip
The above Lemma shows that the operator $C$ is a singular perturbation
of $-d^2/dx^2$ by a delta function at the point $0$ studied
also in [1]. We prove now the stated result

\proclaim Theorem {5.3}. The minimal q.d.s. constructed from
the above operators $G$ and $L_\ell$ is conservative whenever
$g\in H^1((0,+\infty);\IC)$.

\noindent{\bf Proof.}\quad We check assumption {\bf C} for the operator
$C=-2G$ in order to apply Theorem 4.4. All vectors $u\in D(G)$  belong
to $H^2((0,+\infty);\IC)$; hence $L_1u$ belongs to $H^1((0,+\infty);\IC)$.
Moreover, by Lemma 5.2, we have
$$\eqalign{
2\re\left\langle Cu,Gu\right\rangle
+\sum_{\ell =1}^2 \left\Vert C^{1/2}L_\ell u\right\Vert^2
&= - \left\Vert u''\right\Vert^2
+\sum_{\ell=1}^2 \left(\left\Vert L_1(L_\ell u)\right\Vert^2
+ \left\Vert L_2(L_\ell u)\right\Vert^2\right)\cr
&= {|u(0)|^2\over 2\alpha}\left(\left\Vert g'\right\Vert^2
+{\theta^2}\left\Vert g\right\Vert^2
+{|g(0)|^2\over 2\alpha}\right).\cr}
$$
Thus, when $g=0$, assumption {\bf C} obviously holds. If $g\not=0$
then Lemma 5.2 yields the inequality
$$
(2\alpha)^{-1}|u(0)|^2
= \left\Vert g\right\Vert^{-2}\left\Vert L_2 u\right\Vert^2
\le \left\Vert g\right\Vert^{-2}\left\Vert C^{1/2}u\right\Vert^2.
$$
Hence we have
$$
2\re\left\langle Cu,Gu\right\rangle
+\sum_{\ell =1}^2 \left\Vert C^{1/2}L_\ell u\right\Vert^2
\le \left({\left\Vert g'\right\Vert^2\over \left\Vert g\right\Vert^2}
+{\theta^2}+{|g(0)|^2\over 2\alpha\left\Vert g\right\Vert^2}\right)
\left\Vert C^{1/2}u\right\Vert^2.
$$
Therefore assumption {\bf C} holds. The proof is complete letting
$\Phi=C$ and applying Theorem 4.4.\quad\qed

\medskip
\noindent{\bf 5.3 Non closable forms.}
\par\noindent
We study a minimal q.d.s. constructed from operators
$G$, $L_\ell$ so singular that the quadratic form
$u\to-2\re\left\langle u,Gu\right\rangle$ with domain $D(G)$ is not
closable. This problem also can not be solved applying the tools
developed in [9], [11].

Let us consider the contraction semigroup in $h=L^2(0,+\infty)$
$$
P(t)u(x) = u(x+t)
$$
with infinitesimal generator $G$ given by
$$
D(G) = H^1(0,+\infty),\qquad\quad Gu = u'
$$
Let $L_\ell=0$ for $\ell\ge 2$ and let $L_1$ be the operator in $h$
$$
D(L) = H^1(0,+\infty),\qquad\quad Lu = u(0) g
$$
where $g\in h$ and $\Vert g\Vert=1$. Clearly condition {\bf A} holds.
Let $C$ be the self-adjoint operator in $h$
$$
D(C) = \left\{ u\in H^2(0,+\infty)\mid u'(0)=u(0)\right\},
\qquad Cu = - 2 u'',
$$
Applying Lemma 5.2 we can prove that:
\item{(a)} the domain of $C^{1/2}$ contains $H^1(0,+\infty)=D(G)$,
\item{(b)} for all $u\in H^1(0,+\infty)=D(G)$ we have
$$
-2\re\langle u, Gu\rangle
=  |u(0)|^2 \le \Vert C^{1/2}u\Vert^2
=2 \left(\Vert u'\Vert^2 + |u(0)|^2\right),
$$
\item{(c)} for all $u\in D(G^2)=H^2(0,+\infty)$ and
$g\in H^1((0,+\infty);\IC)$ we have
$$\eqalign{
2\re\left\langle Cu,Gu\right\rangle
+\left\Vert C^{1/2}Lu\right\Vert^2
&= -2\left\langle u'',u'\right\rangle -2\left\langle u',u''\right\rangle
+ 2|u(0)|^2 \left\Vert C^{1/2}g\right\Vert^2\cr
&= 2|u(0)|^2\left( 1+ |g(0)|^2 + \Vert g'\Vert^2\right)\cr
&\le 2\left( 1 + |g(0)|^2 + \Vert g'\Vert^2\right)
\left\Vert C^{1/2}u\right\Vert^2.\cr}
$$
Therefore conditions {\bf A} and {\bf C} hold whenever
$g\in H^1((0,+\infty);\IC)$.

By Theorem 4.3, in order to show that the minimal q.d.s. constructed
>from the above operators $G$ and $L$ is conservative it suffices to
check the inequality
$$
\left\langle u,F_n u\right\rangle
\le 2\left(\left\Vert u'\right\Vert^2 + |u(0)|^2\right)
= \left\langle u,Cu\right\rangle
$$
for $u\in D(C)$, $n\ge 1$ where $F_n$ is the unique bounded
extension of $|nLR(n;G)|^2$. A straightforward computation yields
$$
\left\Vert nLR(n;G)u\right\Vert^2
=\left|n\int_0^\infty\hbox{e}^{-nt}u(t) dt\right|^2
$$
Integrating by parts for $u\in D(C)$ and using the Schwarz
inequality we have
$$
\left|n\int_0^\infty\hbox{e}^{-nt}u(t) dt\right|^2
= \left|u(0) +\int_0^\infty \hbox{e}^{-nt}u'(t) dt\right|^2
\le 2 |u(0)|^2 + {1\over n}\left\Vert u'\right\Vert^2.
$$
Therefore the hypotheses of Theorem 4.3 are satisfied and
the minimal q.d.s. is conservative.

It is worth mentioning here that the restriction of this
q.d.s. to the abelian algebra of multiplication operators
by a bounded function coincides with the infinitesimal generator
of a classical stochastic process that can be described as follows:
a point moves on $]0,+\infty[$ towards $0$ with constant speed,
when it reaches $0$ it jumps back on an interval $(a,b)$ with
probability $\left\Vert g 1_{(a,b)}\right\Vert^2$. Journ\'e
showed in [18] that the above minimal q.d.s. is conservative
for every $g\in h$.

\vskip 0.6truecm
\noindent{\bf Comment added in proof.} {\cmrnine
The use of the ``reference'' operator $C$ in (4.2) dominating $\Phi$
was done independently also by A.S.~Holevo in [24].
We arrive to this observation from the resolvent
analyses and we find this assumption quite relevant from the
physical point of view because it allows to deal with
noncomparable operators $H$ and $\Phi$. However in [24]
additional hypotheses were used: (1) the form
$u\to-2\Re e\langle u,Gu\rangle$ is closed, (2) the operators $L_\ell$
are closed, (3) $\Vert Hu\Vert \le\Vert Cu\Vert$
for $u$ in a common core for $G$, $G^*$ and $C$.
The last assumption means essentially that $H$ is dominated by $C$.
On the other hand tangible interpretations of assumptions (1)
and (2), which are not fulfilled in Example 5.2 and 5.3, is unclear. }

\vskip 0.6truecm
\centerline{\bf References}

\item{[1]} S.~Albeverio, F.~Gesztesy, R.~Hoegh-Krohn, H.~Holden:
{\sl Solvable Models in Quantum Mechanics}, Springer,
Berlin-Heidelberg-New York, (1988).

\item{[2]}{L. Accardi, R. Alicki, A. Frigerio, Y.G. Lu}: An invitation
to the weak coupling and low density limits. {\sl Quantum Probability and
related topics} {\bf VI}, (1991), 3--61.

\item{[3]} R.~Alicki, A.~Frigerio: Scattering theory for quantum
dynamical semigroups II. {\sl Ann. Inst. Henri Poincar\'e}
{\bf XXXVIII}, n. 2, (1983), 187-197.

\item{[4]} R.~Alicki, K.~Lendi: {\sl Quantum dynamical semigroups
and applications}.  {\sl Lect. Notes Phys.} {\bf 286},
Springer Verlag, Berlin Heidelberg, New York (1987).

\item{[5]} B.V.R.~Bhat, F.~Fagnola, K.B.~Sinha: On quantum
extensions of semigroups of brownian motions on an half-line.
{\sl Russian J. Math. Phys.} {\bf 4}, (1996),  13--28.

\item{[6]}
B.V.R.~Bhat, K.B.~Sinha: Examples of unbounded
generators leading to non-conservative minimal semigroups.
{\sl Quantum Probability and Related Topics} {\bf IX} (1994), 89--103.

\item{[7]} A.M.~Chebotarev: Sufficient conditions for conservativity
of dynamical semigroups. {\sl Theor. Math. Phys.} {\bf 80}, 2 (1989).

\item{[8]} A.M.~Chebotarev: {\sl The theory of conservative
dynamical semigroups and its applications.} Preprint MIEM n.1.
March 1990.

\item{[9]} A.M.~Chebotarev: Sufficient conditions of the conservativism
of a minimal dynamical semigroup. {\sl Math. Notes} {\bf 52}, (1993),
1067--1077.

\item{[10]} A.M.~Chebotarev: Application of quantum probability to
classical stochastics. Univ. degli Studi di Roma "Tor Vergata",
Centro V.~Volterra, March 1996, Preprint N 246.

\item{[11]} A.M.~Chebotarev, F.~Fagnola: Sufficient Conditions for
Conservativity of Qu\-antum Dynamical Semigroups. {\sl J. Funct. Anal.}
{\bf 118} (1993), 131--153.

\item{[12]} E.B.~Davies: Quantum dynamical semigroups and the neutron
diffusion equation. {\sl Rep. Math. Phys.} {\bf 11} (1977), 169--188.

\item{[13]} F.~Fagnola: Chebotarev's sufficient conditions for
conservativity of quantum dynamical semigroups. {\sl Quantum
Probability and Related Topics} {\bf VIII} (1993) 123--142.

\item{[14]} F.~Fagnola: Characterization of isometric and unitary
weakly differentiable cocycles in Fock space.
Preprint UTM n.358, Trento, October 1991. {\sl Quantum
Probability and Related Topics} {\bf VIII}, (1993), 143--164.

\item{[15]} F.~Fagnola: Diffusion processes in Fock space.
{\sl Quantum Probability and Related Topics} {\bf IX} (1994), 189--214.

\item{[16]} A.S.~Holevo: On the Structure of Covariant Dynamical
Semigroups. {\sl J. Funct. Anal.} {\bf 131} (1995), 255--278.

\item{[17]} K.~Ichihara: Explosion problems for symmetric diffusion
processes. {\sl in} ``Sto\-chas\-tic Processes and Their Applications''
(K.~It\^o and T.~Hida Eds.) {\sl Lect. Notes Math.} {\bf 1203}
(1986), 75--89.

\item{[18]} J.-L.~Journ\'e: Structure des cocycles markoviens sur
l'espace de Fock. {\sl Prob\-ab. Th. Rel. Fields}
{\bf 75}, (1987), 291--316.

\item{[19]} T.~Kato: {\sl Perturbation Theory for Linear Operators.}
Springer. 1966.

\item{[20]}
P.-A.~Meyer: {\sl Quantum Probability for Probabilists}.  {\sl Lect. Notes
Math.} {\bf 1538}. Springer Verlag, Berlin Heidelberg, New York 1993.

\item{[21]} K.R. Parthasarathy: {\sl An Introduction to Quantum
Stochastic Calculus.} Mo\-no\-graphs in Mathematics, Birkh\"auser,
Basel 1992.

\item{[22]} M. Reed, B. Simon: {\sl Methods of Modern Mathematical
Physics.} Vol. I: Functional Analysis. Academic Press.
New York and London 1972.

\item{[23]} Kalyan B.~Sinha: Quantum Dynamical Semigroups. {\sl in}
``Operator Theory: Advances and Applications'', Vol. 70, Birkh\"auser
Verlag, Basel, 1994.

\item{[24]} A.S.~Holevo: Stochastic differential equations in Hilbert
space and quantum Markovian evolution. In: S. Watanabe, M. Fukujima,
Yu. V. Prokhorov, A.N. Shiryaev (eds.) {\sl Proceedings of the VII
Japan-Russia Symposium}, p. 122--131, World Scientific 1996.

\bye